\documentstyle[11pt]{article}
\newtheorem {th}{Theorem}[section]
\newtheorem {pr}[th]{Proposition}

\newtheorem {cor}[th]{Corollary}
\newtheorem {defn}[th]{Definition}

\def\pa{{\tt Path}}
\def\Cox{\hfill \Box}

\def\sf{\sigma\mbox{-field}}

\def\ee{\epsilon}
\def\E{{\bf{E}}}
\def\P{{\bf{P}}}
\def\N{\hbox{I\kern-.2em\hbox{N}}}
\def\R{\hbox{I\kern-.2em\hbox{R}}}

\def\Z{{\bf{Z}}}

\def\B{{\cal{B}}}
\def\|{\, | \, }
\def\v0{{\bf 0}}
\def\one{{\bf 1}}
\def\0{\hat{0}}
\def\1{\hat{1}}

\def\Gt{{\tilde{G}}}

\def\Kt{{\tilde{K}}}

\begin{document}
\begin{center}
{\large \bf MARTIN CAPACITY FOR MARKOV CHAINS } \\
\end{center}
\vspace{2ex}
\begin{center}
{\sc Itai Benjamini}\footnote{
 Research partially supported by the U. S.
Army Research Office through the Mathematical Sciences Institute
of Cornell University.}, \, 
{\sc Robin Pemantle} \footnote{Research supported in part by 
National Science Foundation grant \# DMS 9300191, by a Sloan Foundation
Fellowship, and by a Presidential Faculty Fellowship.}
\,  and \, 
{\sc Yuval Peres} \footnote {Research partially supported by NSF grant
\# DMS-9404391 and a Junior Faculty Fellowship from the Regents of
the University of California\,.}
\end{center}
\begin{center}
\it Cornell University, University of Wisconsin and University of California 
\rm
\end{center}
\vspace{5ex}
\begin{center}
{\bf Abstract}
\end{center}
The probability that a transient Markov chain, or a Brownian path, 
will ever visit a given set $\Lambda$, is classically estimated using
the capacity of $\Lambda$ with respect to the Green kernel $G(x,y)$.
We show that  replacing the Green kernel by the Martin kernel 
 $G(x,y)/G(0,y)$ yields
improved estimates, which are exact up to a factor of 2. These estimates
are applied to random walks on lattices, and also to explain a connection
found by R. Lyons between capacity and percolation on trees.
\vfill



 
\newpage

\section{Introduction}

Kakutani (1944) discovered that a compact set $\Lambda \subseteq \R^d$ is
hit with positive probability by a $d$-dimensional Brownian motion
$(d \geq 3)$ if and only if $\Lambda$ has positive Newtonian 
capacity.  A more quantitative relation holds between this
probability and capacity.  Under additional assumptions on
the set $\Lambda$ it is well known that the hitting probability
is estimated by capacity up to a constant factor (often unspecified).  
Our main result, Theorem~\ref{th2.1}, estimates this probability
(for general Markov chains) by a different capacity up to a factor of 2.  
We first state this estimate for Brownian motion.

\begin{pr} \label{pr1.1}
Let $\{ B_d (t) \}$ denote standard $d$-dimensional Brownian 
motion with $B_d (0) = \v0$ and $d \geq 3$.  Let $\Lambda$ be any closed set
 in $\R^d$.  Then
\begin{equation} \label{eq1.1}
{1 \over 2} {\rm Cap}_K (\Lambda) \leq \P [ \exists t > 0 : B_d (t) \in 
    \Lambda ] \leq {\rm Cap}_K (\Lambda)
\end{equation}
where 
$$K (x,y) = {||y||^{d-2} \over ||x-y||^{d-2}} $$
for $x \neq y$ in $\R^d$, and $K(x,x)=\infty \, $.
  Here $||x-y||$ is the Euclidean distance and
$${\rm Cap}_K (\Lambda) = \left [ \inf_{\mu (\Lambda) = 1} \int_\Lambda
    \int_\Lambda K(x,y) \, d\mu(x) \, d\mu (y) \right ]^{-1} .$$
\end{pr}

\noindent{\bf Remarks:} \\ 
1.\ More detailed definitions will be given in the next section.\\
2.\ The constants $1/2$ and $1$ in~(\ref{eq1.1}) are sharp (see Section 4).

Note that while the Green kernel $G(x,y) = ||x-y||^{2-d}$,
and hence the corresponding
capacity, are translation invariant, the hitting probability of a set
$\Lambda$ by standard $d$-dimensional Brownian motion is not translation 
invariant, but is invariant under scaling. This scale-invariance is shared
by the Martin kernel $K(x,y) = G(x,y) / G(0,y)$. 

The rest of this paper is organized as follows.  Section~2 states
and proves the connection between the probability of a Markov chain 
hitting a set and the  Martin capacity of the set.  Section~3 
gives several examples, including a relation between simple random walk
in three dimensions and the time-space chain arising from simple
random walk in the plane.  
The ratio of $2$ between the two sides in the estimate (\ref{eq1.1})
may remind the reader of a theorem of Lyons (1992) that gives
a precise relation between capacity and independent percolation on trees.
In SectionT~4 we show how recognizing a ``hidden'' Markov chain in the
percolation setting leads to a very short proof of this theorem.
In Section~5 we give the  proof of Proposition~\ref{pr1.1} 
concerning Brownian motion.  
Section~6 discusses motivations and extensions.

\section{Main result}

First we recall some potential theory notions.

\begin{defn} \label{def1}
Let $\Lambda$ be a set and $\B$ a $\sf$ of subsets of $\Lambda$.  
Given a measurable function $F : \Lambda \times \Lambda \rightarrow
[0,\infty]$ and a finite measure $\mu$ on $(\Lambda , \B)$, the
$F$-energy of $\mu$ is 
$$I_F (\mu) = \int_\Lambda \int_\Lambda F(x,y) \, d\mu (x) \, 
   d\mu (y) .$$
The capacity of $\Lambda$ in the kernel $F$ is
$${\rm Cap}_F (\Lambda) = \left [ \inf_\mu I_F (\mu) \right ]^{-1}$$
where the infimum is over probability measures $\mu$ on $(\Lambda , \B)$ 
and by convention, $\infty^{-1} = 0$.  
\end{defn}

If $\Lambda$ is contained in Euclidean space, we always take $\B$ to
be the Borel $\sf$; if $\Lambda$ is countable, we take $\B$
to be the $\sf$ of all subsets.  When $\Lambda$ is countable we also
define the {\em asymptotic capacity} of $\Lambda$ in the kernel $F$:
\begin{equation} \label{eq2.1}
{\rm Cap}_F^{(\infty)} (\Lambda) = \inf_{\{\Lambda_0 \mbox{ finite}\}}
 {\rm Cap}_F (\Lambda
    \setminus \Lambda_0) .
\end{equation}

Let $\{ p(x,y) : x,y \in Y \}$ be transition probabilities on the countable
set $Y$, i.e.\ $\sum_y p(x,y) = 1$ for every $x \in Y$.  Let $\rho \in
Y$ be a distinguished starting state and let $\{ X_n : n \geq 0 \}$
be a Markov chain with $\P [ X_{n+1} = y \| X_n = x ] = p(x,y)$.

Define the Green function 
$$G(x,y) = \sum_{n=0}^\infty p^{(n)} (x,y) = \sum_{n=0}^\infty 
   \P_x [X_n = y] $$
where $p^{(n)} (x,y)$ are the $n$-step transition probabilities and
$\P_x$ is the law of the chain \newline
$\{ X_n : n \geq 0 \}$ when $X_0 = x$.
We want to estimate the
probability that a sample path $\{ X_n \}$ visits a set $\Lambda \subseteq Y$.
We assume that the Markov chain $\{ X_n \}$ is transient; in fact, it suffices
to assume that $G(x,y) < \infty$ for all $x,y \in \Lambda$. 

\begin{th} \label{th2.1}
Let $\{ X_n \}$ be a transient Markov chain on the countable state space 
$Y$ with initial state $\rho$ and transition probabilities $p(x,y)$.
For any subset $\Lambda$ of $Y$ we have
\begin{equation} \label{eq2.2}
{1 \over 2} {\rm Cap}_K (\Lambda) \leq \P_\rho [\exists n \geq 0 : X_n \in \Lambda]
   \leq {\rm Cap}_K (\Lambda)
\end{equation}
and
\begin{equation} \label{eq2.3}
{1 \over 2} {\rm Cap}_K^{(\infty)} (\Lambda) \leq \P_\rho [
   X_n \in \Lambda \mbox{ infinitely often }] \leq {\rm Cap}_K^{(\infty)} (\Lambda)
\end{equation}
where $K$ is the Martin kernel 
\begin{equation} \label{eq2.4}
K(x,y) = {G(x,y) \over G(\rho , y)}
\end{equation}
defined using the initial state $\rho$.
\end{th}
\noindent{\bf Remarks:}
\begin{enumerate}
\item The Martin kernel $K(x,y)$ can obviously be replaced by the symmetric
  kernel \newline $\frac{1}{2} (K(x,y)+K(y,x))$ without affecting the
   energy of measures or the capacity of sets.
\item If the Markov chain starts according to an initial measure
 $\pi$ on the state space, rather than from a fixed initial state,
 the theorem may be applied by adding an abstract initial state $\rho$
 with transition probabilities $p(\rho,y) = \pi(y)$ for $y \in Y$.
\end{enumerate}

{\noindent \sc Proof:}  $(i)$ The right hand inequality in~(\ref{eq2.2})
follows from an entrance time decomposition.  Let $\tau$ be the 
first hitting time of $\Lambda$ and let $\nu$ be the (possibly
defective) hitting measure
$\nu (x) = \P_\rho [X_\tau = x]$ for $x \in \Lambda$.  Then 
\begin{equation} \label{eq2.5}
\nu (\Lambda) = \P [ \exists n \geq 0 : X_n \in \Lambda] \, .
\end{equation}
Now for all $y \in \Lambda$ :
$$\int G(x,y) \, d\nu (x) = \sum_{x \in \Lambda} \P_\rho [ X_\tau = x]
    G(x,y) = G(\rho , y) .$$
Thus $\int K(x,y) \, d\nu (x) = 1$ for every $y \in \Lambda$.  
Consequently 
$$I_F \left ( {\nu \over \nu (\Lambda)} \right ) = \nu (\Lambda)^{-2}
   I_F (\nu) = \nu (\Lambda)^{-1},$$
so that ${\rm Cap}_K (\Lambda) \geq \nu (\Lambda)$.  By~(\ref{eq2.5}), this
proves half of~(\ref{eq2.2}).

To establish the left hand inequality in~(\ref{eq2.2}) we use the
second moment method.  Given a probability measure $\mu$ on 
$\Lambda$, consider the random variable
$$Z = \int_\Lambda G(\rho , y)^{-1} \sum_{n=0}^\infty \one_{\{X_n = y\}}
   \, d\mu (y) .$$
By Tonelli and the definition of $G$, 
\begin{equation} \label{eq2.6}
\E_\rho Z = 1.
\end{equation}
Now we bound the second moment:
\begin{eqnarray*}
\E_\rho Z^2 & = & \E_\rho \int_\Lambda \int_\Lambda G(\rho,y)^{-1} 
   G(\rho , x)^{-1} \sum_{m,n = 0}^\infty \one_{\{X_m = x , X_n = y\}}
   \, d\mu (x) \, d\mu (y) \\[2ex]
& \leq & 2 \E_\rho \int_\Lambda \int_\Lambda G(\rho,y)^{-1} 
   G(\rho , x)^{-1} \sum_{0 \leq m \leq n < \infty} \one_{\{X_m = x , 
   X_n = y\}} \, d\mu (x) \, d\mu (y) .
\end{eqnarray*}
For each $m$ we have
$$\E_\rho \sum_{n=m}^\infty \one_{\{X_m = x, X_n = y\}} \;=\; \P_\rho [X_m = x]
   G(x,y) .$$
Summing this over all $m \geq 0$ yields $G(\rho , x) G(x,y)$, and therefore
$$\E_\rho Z^2 \leq 2 \int_\Lambda \int_\Lambda G(\rho,y)^{-1}
   G(x,y) \, d\mu (x) \, d\mu (y) = 2 I_K (\mu) .$$
By Cauchy-Schwarz and~(\ref{eq2.6}), 
$$\P_\rho [ \exists n \geq 0: X_n \in \Lambda] \geq \P_\rho [ Z > 0]
   \geq {(\E_\rho Z)^2 \over \E_\rho Z^2} \geq {1 \over 2 I_K (\mu)} .$$
Since the left hand side does not depend on $\mu$, we conclude that
$$\P_\rho [ \exists n \geq 0 : X_n \in \Lambda] \geq {1 \over 2} {\rm Cap}_K
   (\Lambda) $$
as claimed.   

To infer~(\ref{eq2.3}) from~(\ref{eq2.2}) observe that since
$\{ X_n \}$ is a transient chain, almost surely every state is visited
only finitely often and therefore 
$$ \left \{ X_n \in \Lambda ~\mbox{  infinitely often } \right \} = 
   \bigcap_{\Lambda_0 \mbox{ finite}} \left \{ \exists n \geq 0 :
   X_n \in \Lambda \setminus \Lambda_0 \right \} ~\mbox{ a.s.}$$
Applying~(\ref{eq2.2}) and the definition~(\ref{eq2.1}) of asymptotic
capacity yields~(\ref{eq2.3}).   $\Cox$

\section{Corollaries and examples}

This section is devoted to deriving some consequences
of Theorem~\ref{th2.1}.  The first involves a widely applicable 
equivalence relation between distributions of random sets. 

\noindent{\bf Definition}: 
 Say that two random subsets 
$W_1$ and $W_2$ of a countable space are \newline {\em intersection-equivalent}
 (or more precisely,
that their {\em laws} are intersection-equivalent) if there exist
positive finite constants $C_1$ and $C_2$, such that for every subset 
$A$ of the space, 
$$ C_1 \leq {\P [W_1 \cap A \neq \emptyset] \over
             \P [W_2 \cap A \neq \emptyset]} \leq C_2 \, .
$$
 It is easy to see that if $W_1$ and $W_2$
are intersection-equivalent then \newline
$C_1 \leq \P [ |W_1 \cap A| = \infty] / 
          \P [ |W_2 \cap A| = \infty] \leq C_2$ for all sets $A$, with the 
  same constants $C_1$ and $C_2$.
  An immediate corollary 
of Theorem~\ref{th2.1} is the following, one instance of which
is given in Corollary~\ref{cor2.6}.
\begin{cor} \label{cor I-E}
Suppose the Green's functions for two Markov chains on the
same state space {\rm (}with the same initial state{\rm )}
 are bounded by constant multiples of each other.
{\rm (}It suffices that this bounded ratio property holds 
 for the corresponding Martin kernels
 $K(x,y)$ or for their symmetrizations $K(x,y)+K(y,x)$.{\rm )}
Then the ranges of the two chains are intersection-equivalent. 
\end{cor}

Lamperti (1963) gave an alternative criterion for $\{ X_n \}$ to
visit the set $\Lambda$ infinitely often.  Fix $b > 1$.  With the notations 
of Theorem~\ref{th2.1}, denote $Y(n) = \{ x \in Y : b^{-n-1} < 
G(\rho , x) \leq b^{-n} \}$.

\begin{cor}[Lamperti's Wiener Test] \label{cor2.2}
Assume that the set $\{ x \in Y : G(\rho , x) > 1 \}$ is finite. Also, 
assume that there exists a constant $C$ such that
for all sufficiently large $m$ and $n$ we have
\begin{equation} \label{eq2.7}
G(x,y) < C b^{-(m+n)} 
\end{equation}
for all $x \in Y(m)$ and $y \in Y(m+n).$  Then
\begin{equation} \label{eq2.8}
\P [X_n \in \Lambda \mbox{ infinitely often}] > 0 \, \Longleftrightarrow \,
   \sum_{n=1}^\infty b^{-n} {\rm Cap}_G(\Lambda \cap Y(n)) = \infty \, .
\end{equation}
\end{cor}

\noindent{\sc Sketch of proof:} Clearly 
$\sum_{n=1}^\infty b^{-n} {\rm Cap}_G(\Lambda \cap Y(n)) = \infty$ if
and only if $\sum_n {\rm Cap}_K (\Lambda \cap Y(n)) = \infty$.  The
equivalence (\ref{eq2.8}) then follows from a version of the Borel-Cantelli
lemma proved in Lamperti's paper (a better proof is in
Kochen and Stone (1964)). 

 Lamperti's Wiener test is useful in many cases; however the 
condition~(\ref{eq2.7}) excludes some natural transient chains
such as simple random walk on a binary tree.
Next, we deduce from Theorem~\ref{th2.1} a criterion for a
 recurrent Markov chain to visit
its initial state infinitely often within a prescribed time set.
\begin{cor} \label{cor2.3}
Let $\{ X_n \}$ be a recurrent Markov chain on the countable state
space $Y$, with initial state $X_0 = \rho$ and transition probabilities
$p(x,y)$.  For nonnegative integers $m \leq n$ denote
$$\Gt (m,n) = \P [X_n = \rho \| X_m = \rho] = p^{(n-m)} (\rho , \rho)$$
and 
$$\Kt (m,n) = {\Gt (m,n) \over \Gt (0,n)} .$$
Then for any set of times $A \subseteq \Z^+$:
\begin{equation} \label{eq2.10}
{1 \over 2} {\rm Cap}_{\Kt} (A) \leq \P [\exists n \in A : X_n = \rho]
   \leq {\rm Cap}_{\Kt} (A)
\end{equation}
and
\begin{equation} \label{eq2.11}
{1 \over 2} {\rm Cap}_\Kt^{(\infty)} (A) \leq \P [\sum_{n \in A} \one_{\{X_n = \rho\}}
   = \infty] \leq {\rm Cap}_\Kt^{(\infty)} (A) .
\end{equation}
\end{cor}

\noindent{\sc Proof:}  Consider the space-time chain $\{ (X_n , n) :
n \geq 0 \}$ on the state space $Y \times \Z^+$.  This chain is
obviously transient; let $G$ denote its Green function.  Since
$G((\rho , m) , (\rho , n)) = \Gt (m,n)$ for $m \leq n$, applying
Theorem~\ref{th2.1} with $\Lambda = \{\rho\} \times A$ shows 
that~(\ref{eq2.10}) and~(\ref{eq2.11}) follow respectively 
{}from~(\ref{eq2.2}) and~(\ref{eq2.3}).   $\Cox$

\noindent{\bf Example 1:} {\em Random walk on} $\Z$.  The moral of this
example will be that Borel-Cantelli does not always correctly
settle questions about return times of random walks;
similar examples may be found in Ruzsa and Sz\'ekely (1982) 
and Lawler (1991).  

Let $S_n$ be the partial sums
of mean-zero, finite variance, i.i.d.\ integer random variables.  By
the local central limit theorem (cf. Spitzer 1964), 
$$\Gt (0,n) = \P [S_n = 0] \approx c n^{-1/2}$$
provided that the summands $S_n - S_{n-1}$ are aperiodic.  Therefore
\begin{equation} \label{eq2.12}
\P [ \sum_{n \in A} \one_{\{S_n = 0\}} = \infty] > 0 \Leftrightarrow
   {\rm Cap}_F^{(\infty)} (A) > 0,
\end{equation}
with $F(m,n) = (n^{1/2} / (n-m+1)^{1/2}) \one_{\{m \leq n\}}$.
  By the Hewitt-Savage
zero-one law, the event in~(\ref{eq2.12}) must have probability zero
or one.  Consider the special case in which $A$ consists of separated
blocks of integers:
\begin{equation} \label{eq2.14}
A = \bigcup_{n=1}^\infty [2^n , 2^n + L_n] .
\end{equation}
A standard calculation (e.g., with the Wiener test applied to the time-space
chain)  shows that in this case 
$S_n = 0$ for infinitely many $n \in A$ with probability one, if and only if 
$\sum_n L_n^{1/2} 2^{-n/2} = \infty$.  On the other hand, the expected
number of returns $\sum_{n \in A} \P [S_n = 0]$ is infinite if and
only if $\sum_n L_n 2^{-n/2} = \infty$.  Thus an infinite expected number
of returns in a time set does not suffice for almost sure return in
the time set.  When the walk is periodic, i.e.\ 
$$r = \mbox{gcd} \{ n : \P [S_n = 0] > 0 \} > 1,$$
the same criterion holds as long as $A$ is contained in $r \Z^+$.

In some cases, the criterion of Corollary~\ref{cor2.3} can be turned
around and used to estimate asymptotic capacity.  For instance, if
$\{S_n'\}$ is an independent random walk with the same distribution as
$\{ S_n \}$ and $A$ is the random set $A = \{ n : S_n' = 0 \}$,
then the positivity of ${\rm Cap}_F^{(\infty)} (A)$ follows from the 
recurrence of the planar random walk $\{ (S_n , S_n') \}$.  
This  implies that the ``discrete Hausdorff dimension'' of $A$
(in the sense of Barlow and Taylor (1992)) is almost surely $1/2$ ;
detailed estimates of the discrete Hausdorff measure of $A$ were
obtained by Khoshnevisan~(1993).

\noindent{\bf Example 2:} {\em Random walk on} $\Z^2$.  Now we assume 
that $S_n$ are partial sums of aperiodic, mean-zero, finite variance 
i.i.d.\ random variables in
$\Z^2$.  Let $A \subseteq \Z$.  Again, $\P [S_n = 0 \mbox{ for infinitely
many } n \in A]$ is zero or one and it is one if and only if 
${\rm Cap}_F^{(\infty)}(A) > 0$ where 
 $F(m,n) = (n / {(1+n-m)}) \one_{\{m \leq n\}}$. This
follows from the local central limit theorem 
(cf.\ Spitzer 1964) which ensures that 
$$\Gt (0 , n) = \P[S_{n} = 0] \approx cn^{-1} \mbox{ as } n 
   \rightarrow \infty .$$
For instance, if $A$ consists of disjoint blocks 
$$A = \bigcup_n [2^n , 2^n + L_n]$$
then ${\rm Cap}_F^{(\infty)} (A) > 0$ if and only if $\sum_n 2^{-n} L_n / \log L_n
= \infty$.  The expected number of returns to zero is infinite if and only
if $\sum 2^{-n} L_n = \infty$. 

Comparing the kernel $F$ with the Martin kernel for simple random 
walk on $\Z^3$ leads to the next corollary.
\begin{cor} \label{cor2.6}
For $d=2,3 \, $, let $\{ S_n^{(d)} \}$ be a truly d-dimensional
 random walk on the $d$-dimensional
lattice, with  increments of mean zero and finite variance.
Assume that the walks are aperiodic, i.e., the set of positive
integers n for which $\P [S^{(d)}_n = 0] > 0$ has g.c.d. $=1$.
 Then there exist positive finite constants $C_1$ and $C_2$
such that for any set of positive  integers $A$,
\begin{equation} \label{eq2.15}
C_1 \leq \frac{\P [S^{(2)}_{n} = 0 \mbox{ for some } n \in A]}
     {\P [S^{(3)}_n \in {\{0\} \times \{0\} \times A}
 \mbox{ for some } n ]} \leq C_2 \, , 
\end{equation}
where $\{0\} \times \{0\} \times A = \{(0,0,k) : k \in A \}$. Consequently,
\begin{equation} \label{eq2.16}
\P [S^{(2)}_{n} = 0 \mbox{ for infinitely many } n \in A] = 
\P [S^{(3)}_n \in \{0\} \times \{0\} \times A \mbox{ infinitely often}].
\end{equation}
\end{cor}
Note that both sides of (\ref{eq2.16}) take only the values 0 or 1.
Corollary \ref{cor2.6} follows from Corollary \ref{cor I-E}, in conjunction
with Example 2 above and the asymptotics
 $G(0,x) \sim c / |x|$ as 
$|x| \rightarrow \infty$ for the random walk $S_n^{(3)}$ (cf. Spitzer (1964)).
 The Wiener test implies the equality (\ref{eq2.16}) but not the estimate
(\ref{eq2.15}).
  Erd\"os (1961) and McKean (1961) showed that 
for $A = \{ \mbox{primes}\}$
 , the left-hand side of~(\ref{eq2.16}) is 1. The corresponding result for the
right-hand side is in Kochen and Stone (1964).
 To see why Corollary~\ref{cor2.6}
is surprising, observe that the space-time chain $\{ (S_n^{(2)} , n) \}$
travels to infinity faster than $S_n^{(3)}$, yet by Corollary~\ref{cor2.6},
the same subsets of lattice points on the positive $z$-axis are
hit infinitely often by the two processes.  

\noindent{\bf Example 3:} {\em Riesz-type kernels.}
  The analogues of the Riesz kernels
in the discrete setting are the kernels 
$$F_\alpha (x,y) = { ||y||^\alpha \over 1 + ||x-y||^\alpha }$$
on $\Z^d$, where $|| \cdot ||$ is any norm.  We write 
${\rm Cap}_\alpha^{(\infty)}$ for ${\rm Cap}_{F_\alpha}^{(\infty)}$.  By
Theorem~\ref{th2.1}, the asymptotics for the Green function, and
the Hewitt-Savage law, simple random walk on $\Z^d$ visits a set
$\Lambda \subseteq \Z^d$ i.o.\ a.s.\ if and only if ${\rm Cap}_{d-2}^{(\infty)}
(\Lambda) > 0$. More generally, if a random walk $\{S_n\}$
on the $d$-dimensional
lattice has a Green function satisfying
$G(0,x) \sim c |x|^ {\alpha -d}$ as $|x| \rightarrow \infty$, then
Theorem~\ref{th2.1} implies that $S_n \in \Lambda$ for infinitely
many $n$ a.s. iff ${\rm Cap}_{d-\alpha }^{(\infty)} (\Lambda) > 0 $.
These asymptotics for the Green function are known to hold for many 
increment distributions 
in the domain of attraction of an  $\alpha$-stable distribution.
 (cf. Williamson (1968) for some sufficient conditions.)

  Given a set of digits $D \subseteq \{ 0 , 1 , \ldots , 
b-1 \}$ containing zero, consider ``the integer Cantor set'' 
$$\Lambda (D,b) = \{ \sum_{n=0}^N a_n b^n : a_n \in D \, 
  \mbox{\rm for all} \, n\, , \mbox{and} \,  N \geq 0   \} .$$
It may be shown that ${\rm Cap}_\alpha^{(\infty)} (\Lambda (D,b)) > 0$ 
if and only if $|D| \geq b^\alpha$.  This, together with Example 3, 
motivates defining the (discrete) dimension 
of $\Lambda \subseteq \Z^d$ by
\begin{equation} \label{eq2.17}
\dim (\Lambda)= \inf \{ \alpha : {\rm Cap}_\alpha^{(\infty)} (\Lambda) = 0 \} .
\end{equation}
Corollary~8.4 in Barlow and Taylor (1992) shows that this definition 
is equivalent to the definition of discrete Hausdorff dimension  in that paper.

When applying Theorem \ref{th2.1}, it is often useful to know whether for
the Markov chain under consideration, the probability of visiting
a set infinitely often must be either $0$ or $1$.
As remarked before,  random walks on $\Z^d$ (or any abelian group) 
have this property by the Hewitt-Savage zero-one law.
 Easy examples show that this fails for random walk on a 
free group.  More generally, the following ``folklore'' criterion holds.

\begin{pr} \label{pr2.8} 
Let $\mu$ be a probability measure whose support generates a countable 
group $Y$,
and let $\{ S_n \}$ be the random walk with step distribution
$S_n S_{n-1}^{-1} \sim \mu$.  Then the probability \\
$\P [S_n \in \Lambda \mbox{ infinitely often}]$ takes only the
values $0$ and $1$ as $\Lambda$ ranges over subsets of $Y$, if and
only if every bounded $\mu$-harmonic function is constant.  (Recall
that $h : Y \rightarrow \R$ is $\mu$-harmonic if $h(x) = \int_Y
h(yx) \, d\mu (y)$ for all $x \in Y$.)  
\end{pr}

\noindent{\bf Remark:} When all bounded harmonic functions are constant,
one says that the {\em Poisson boundary} of $(Y,\mu)$ is trivial; see
Kaimanovich and Vershik (1982) for background. 

\noindent{\sc Proof:} Given a set $\Lambda \subseteq Y$, the function $h(x)
= \P[S_n x \in \Lambda \mbox{ infinitely often}]$ is bounded and
$\mu$-harmonic.  The Markov property and the martingale convergence
theorem imply that
$$h(S_m) = \P \left [ \{S_k : k \geq 0 \} \mbox{ visits } ~\Lambda 
   \mbox{ i.o. } \| S_1 , S_2 , \ldots , S_m \right ] \; \rightarrow \; 
   \one_{\{S_k \mbox{ visits } {\displaystyle \Lambda} \mbox{ i.o.}\}}$$
as $m \rightarrow \infty$.  Thus if all bounded harmonic functions
are constant, the zero-one law holds.  Conversely, assume the zero-one
law holds and let $h$ be a bounded $\mu$-harmonic function.  For 
$\alpha \in \R$, let $\Lambda_\alpha = \{ y \in Y : h(y) < \alpha \}$.  
If $\P [S_n x \mbox{ visits } \Lambda_\alpha \mbox{ i.o.}] = 0$
then we consider the stopping time $\tau = \min \{n : h(S_n x) 
\geq \alpha \}$ and obtain $h(x) = h(S_0 x) = \E h(S_\tau x) \geq \alpha$.
Similarly, if $\P [S_n x \mbox{ visits } \Lambda_\alpha \mbox{ i.o.}] = 1$  
then $h(x) \leq \alpha$.  Since the support of $\mu$ generates $Y$,
$$\P [S_n \mbox{ visits } \Lambda_\alpha \mbox{ i.o.}] = 0 \Leftrightarrow
  \P [S_n x \mbox{ visits } \Lambda_\alpha \mbox{ i.o.}] = 0 $$
and it follows that $h(x) = h(e)$ for all $x \in Y$.   $\Cox$

\section{Independent percolation on trees}

Theorem \ref{th2.1} yields a short proof of a fundamental result
of R. Lyons concerning  percolation on trees.
This theorem and its variants have been used in the analysis of
a variety of probabilistic processes on trees, including
random walks in a random environment, first-passage percolation
and the Ising model. (See Lyons (1989, 1990, 1992); Lyons
and Pemantle (1992); Benjamini and Peres (1994); Pemantle and Peres (1994).)
Recently, this theorem has also been applied to
 study intersections of sample paths in Euclidean space (cf. Peres (1994)).

\noindent{\bf Notation}: Let $T$ be a finite, rooted tree. 
Vertices of degree one in $T$ (apart from the root $\rho$)
are called {\em leaves}, and the set of leaves is the {\em boundary}
$\, \partial T$ of $T$. The set of edges on the path
connecting the root to a leaf $x$ is denoted $\pa (x)$. 

Independent percolation on $T$ is defined as follows.
To each edge $e$ of $T$, a parameter $p_e$ in $[0,1]$
is attached, and $e$ is removed with probability $1-p_e$,
retained with probability $p_e$, with mutual independence among edges. 
Say that a leaf $x$ {\em survives the percolation} if all of
$\pa (x)$ is retained, and say that the tree boundary 
$\partial T$ survives if
some leaf of $T$ survives.
\begin{picture}(300,350)(-60,-65)
\put(120,270){\line(1,-1){70}}
\put(120,270){\circle*{3}}
\put(120,270){\line(-1,-1){70}}
\put(50,200){\line(0,-1){60}}
\put(50,200){\circle*{3}}
\put(190,200){\line(2,-3){40}}
\put(190,200){\circle*{3}}
\put(190,200){\line(-2,-3){40}}
\put(50,140){\line(1,-2){30}}
\put(50,140){\circle*{3}}
\put(50,140){\line(-1,-2){30}}
\put(150,140){\line(0,-1){60}}
\put(150,140){\circle*{3}}
\put(230,140){\line(1,-2){30}}
\put(230,140){\circle*{3}}
\put(230,140){\line(-1,-2){30}}
\put(20,80){\line(-1,-3){20}}
\put(20,80){\circle*{3}}
\put(20,80){\line(1,-3){20}}
\put(80,80){\line(-1,-3){20}}
\put(80,80){\circle*{3}}
\put(80,80){\line(1,-3){20}}
\put(150,80){\line(0,-1){60}}
\put(150,80){\circle*{3}}
\put(200,80){\line(0,-1){60}}
\put(200,80){\circle*{3}}
\put(260,80){\line(-1,-3){20}}
\put(260,80){\circle*{3}}
\put(260,80){\line(1,-3){20}}
\put(0,20){\circle*{3}}
\put(0,20){\line(-1,-4){10}}
\put(0,20){\line(1,-4){10}}
\put(40,20){\circle*{3}}
\put(40,20){\line(0,-1){40}}
\put(60,20){\circle*{3}}
\put(60,20){\line(-1,-4){10}}
\put(60,20){\line(1,-4){10}}
\put(100,20){\circle*{3}}
\put(100,20){\line(-1,-4){10}}
\put(100,20){\line(1,-4){10}}
\put(150,20){\circle*{3}}
\put(150,20){\line(0,-1){40}}
\put(200,20){\circle*{3}}
\put(200,20){\line(-1,-4){10}}
\put(200,20){\line(1,-4){10}}
\put(240,20){\circle*{3}}
\put(240,20){\line(0,-1){40}}
\put(280,20){\circle*{3}}
\put(280,20){\line(0,-1){40}}
\put(-10,-20){\circle*{3}}
\put(10,-20){\circle*{3}}
\put(40,-20){\circle*{3}}
\put(50,-20){\circle*{3}}
\put(70,-20){\circle*{3}}
\put(90,-20){\circle*{3}}
\put(110,-20){\circle*{3}}
\put(150,-20){\circle*{3}}
\put(190,-20){\circle*{3}}
\put(210,-20){\circle*{3}}
\put(240,-20){\circle*{3}}
\put(280,-20){\circle*{3}}
\put(100,-57){Figure 1:  a tree}
\end{picture}

\begin{th}[Lyons (1992)] \label{Lyons}
With the notation above, define a kernel $F$ on $\partial T$ by \newline
$F(x,y)= \prod \{p_e^{-1} \, : \, e \in \pa (x) \bigcap \pa (y) \}$
for $x \neq y$ and 
$F(x,x)= 2 \prod \{p_e^{-1} \, : \, e \in \pa (x)  \}$.
Then 
$$
 {\rm Cap}_F (\partial T) \leq \P [ \partial T \, \, 
  \mbox{survives the percolation} ]
   \leq 2 {\rm Cap}_F (\partial T)
$$
\end{th}
(The kernel $F$ differs from the kernel used in Lyons (1992)
on the diagonal, but this difference is unimportant in all applications).

\noindent {\sc Proof:}
Embed $T$ in the lower half-plane, with the root at the origin.
The random set of $r \geq 0$ leaves that survive the percolation may be
 enumerated from left to right as $V_1, V_2 , \ldots, V_r$.
The key observation is that 
The random sequence $\rho , V_1 , V_2 , \ldots ,V_r, \Delta,
\Delta, \ldots$ \underline{is a Markov chain} on the state space 
 $\partial T \bigcup \{ \rho, \Delta \} $
 (where $\rho$ is the root and $\Delta$ is a  formal  absorbing cemetery). 

 Indeed, {\em given} that $V_k = x$, all the
edges on $\pa (x)$ are retained, so that
 survival of leafs to the right of $x$ is determined by the
 edges strictly to the right of  $\pa (x)$, and is thus conditionally
 independent of  $V_1, \ldots, V_{k-1}$. This verifies the Markov
 property, so Theorem~\ref{th2.1} may be applied. 

The transition probabilities
 for the Markov chain above are complicated, but it is easy to write down the 
 Green kernel.  Clearly, $G(\rho , y) = \P [ y \mbox{ survives 
 the percolation}] = \prod_{e \in \pa (y)} p_e$.  Also,
 if $x$ is to the left of $y$, then $G(x,y)$ is
 equal to the probability that the range of the Markov chain contains
 $y$ given that it contains $x$, which is just the probability of
 $y$ surviving given that $x$ survives. Therefore
$$G(x,y) = \prod_{e \in \pa (y) \setminus \pa (x)} p_e $$
and hence 
$$K(x,y)={G(x,y) \over G(\rho , y)} 
  = \prod_{e \in \pa (x) \cap \pa (y)} p_e^{-1}\, .
$$
Thus $K(x,y)+K(y,x)=F(x,y)$ for all $x,y \in \partial T$,
and Lyons' Theorem follows from Theorem~\ref{th2.1}.   $\Cox$

\noindent{\bf Remark:} The same method of recognizing
a ``hidden'' Markov chain may be used to prove
more general results on random labeling
of trees due to Evans (1992) and  Lyons (1992).

\section{Martin capacity and Brownian motion}

\noindent{\sc Proof of Proposition}~\ref{pr1.1}: To bound from above the
probability of ever hitting
$\Lambda$, consider the stopping time
$\tau = \min \{ t > 0 : B_d (t) \in \Lambda \}$.  The distribution of
$B_d (\tau)$ on the event $\tau < \infty$ is a possibly defective 
distribution $\nu$ satisfying 
\begin{equation} \label{eq3.1}
\nu (\Lambda) = \P [\tau < \infty] = \P [\exists t > 0 : B_d (t) \in
   \Lambda] .
\end{equation}
Now recall the standard formula, valid when $0 < \ee < ||y|| \, $: 
\begin{equation} \label{eq3.2}
\P [\exists t > 0 : ||B_d (t) - y|| < \ee] = {\ee^{d-2} \over ||y||^{d-2}} \, .
\end{equation}
By a first entrance decomposition, the probability in~(\ref{eq3.2}) is
at least
$$
   \P [ ||B_d (\tau) - y|| > \ee \mbox{ and } \exists t > \tau :
   ||B_d (t) - y|| < \ee] \; = \;
   \int_{x : ||x-y|| > \ee} {\ee^{d-2} \over ||x-y||^{d-2}} \, d\nu (x) \, .$$
Dividing by $\ee^{d-2}$ throughout and letting $\ee \rightarrow 0 $,  we obtain
$$\int_\Lambda {d\nu (x) \over ||x-y||^{d-2}} \leq {1 \over ||y||^{d-2}}\, \,  ,$$
i.e. $\int_\Lambda K(x,y) \, d\nu (x) \leq 1$ for all $y \in \Lambda$.
Therefore $I_K (\nu) \leq \nu (\Lambda)$ and thus 
$${\rm Cap}_K (\Lambda) \geq [I_K (\nu / \nu (\Lambda))]^{-1} \geq \nu (\Lambda) \, ,$$
which by~(\ref{eq3.1}) yields the upper bound on the  probability
of hitting $\Lambda$.  

To obtain a lower bound for this probability, a second moment estimate
is used.  It is easily seen that the Martin capacity of $\Lambda$ is
 the supremum of the capacities of its compact subsets, so we may assume that
 $\Lambda$ itself is compact.
 For $\ee > 0$ and $y \in \R^d$ let $D(y,\ee)$ denote the 
Euclidean ball of radius $\ee$ about $y$ and let $h_\ee (||y||)$ denote
the probability that a standard Brownian path will hit this ball:
\begin{equation} \label{eq3.3}
h_\ee (r) =  \left\{ \begin{array}{ll}
                     ( {\ee / r}  )^{d-2} & \mbox{if $r > \ee$} \\
                      1                   & \mbox{otherwise \, \, .}
                       \end{array} \right.
\end{equation}
Given a probability measure $\mu$ on $\Lambda$, and $\ee > 0$, consider
the random variable 
$$Z_\ee = \int_\Lambda \one_{\{\exists t>0 : B_d(t) \in D(y,\ee)\}} 
   h_\ee (||y||)^{-1} \, d\mu (y)\, .$$  
Clearly $\E Z_\ee = 1$.  We compute the second moment of $Z_\ee$ 
in order to apply Cauchy-Schwarz as in the proof of Theorem~\ref{th2.1}.  

By symmetry, 
\newpage
\begin{eqnarray} \label{eqbrown}
\E Z_\ee^2 & = & 2 \E \int_\Lambda \int_\Lambda \one_{\{\exists t > 0 :
   B_d (t) \in D(x,\ee) \mbox{ and } \exists s > t : B_d (s) \in 
   D(y,\ee)\}} {d\mu(x) d\mu(y) \over h_\ee (||x||) h_\ee (||y||)}
                       \nonumber \\[2ex]
& \leq & 2 \E \int_\Lambda \int_\Lambda \one_{\{\exists t > 0 :
   B_d (t) \in D(x,\ee)\}} {h_\ee (||y-x|| - \ee) \over h_\ee (||x||) 
   h_\ee (||y||)}  \, d\mu (x) \, d\mu (y) \nonumber \\[2ex]
& = & 2 \int_\Lambda \int_\Lambda {h_\ee (||y - x|| - \ee) \over
   h_\ee(||y||)} \, d\mu (x) \, d\mu (y) \, .
\end{eqnarray}
The last integrand is bounded by  1 if $||y|| \leq \ee$. On the other hand, if
$||y|| > \ee $ and $||y-x|| \leq 2 \ee $ then 
 $ h_\ee (||y-x||-\ee) = 1 \leq 2^{d-2} h_\ee (||y-x||) \, $, so that the
integrand on the right-hand side of (\ref{eqbrown}) is at most 
 $2^{d-2} K(x,y) \, $. Thus
\begin{eqnarray} \label{eqbrown2}
&& \E Z_\ee^2 \leq   2 \mu (D(0,\ee)) +
2^{d-1} \int \int \one_{\{||y-x|| \leq 2\ee\}} K(x,y) \, d\mu
   (x) \, d\mu (y)  \\
&& + 2 \int \int \one_{\{||y-x|| > 2\ee\}} \left ( {||y|| \over 
   ||y-x|| - \ee} \right )^{d-2} \, d\mu (x) \, d\mu (y) \, . \nonumber
\end{eqnarray}
Since the kernel is infinite on the diagonal, any measure with
finite energy must have no atoms. Restricting attention to such measures
$\mu$, we see that the  first two summands in (\ref{eqbrown2}) 
drop out as $\ee \rightarrow 0$ (by 
dominated convergence) . This leaves
\begin{equation} \label{eq 20b}
\lim_{\ee \downarrow 0} \E Z_\ee^2 \leq 2 I_K (\mu) \, .
\end{equation}
 Clearly the hitting probability
$\P [\exists t>0 , y \in \Lambda : B_d (t) \in D(y,\ee)]$
is at least 
$$ \P [ Z_\ee > 0] \geq {(\E Z_\ee)^2 \over \E Z_\ee^2} =
   (\E Z_\ee^2 )^{-1} \, .$$
Transience of Brownian motion implies that if the Brownian path 
visits every $\ee$-neighborhood of the compact set $\Lambda$ then it 
almost surely intersects $\Lambda$ itself.  Therefore, by~(\ref{eq 20b}):
$$\P [ \exists t > 0 : B_d (t) \in \Lambda] \geq
   \lim_{\ee \downarrow 0} (\E Z_\ee^2)^{-1} \geq {1 \over 2I_K (\mu)} .$$
Since this is true for all probability measures $\mu$ on $\Lambda$,
we get the desired conclusion:
\begin{equation} \label{eq3.4}
\P [ \exists t > 0 : B_d (t) \in \Lambda] \geq {1 \over 2} {\rm Cap}_K (\Lambda) .
\end{equation}
$\Cox$

\noindent{\bf Remark:}  The right--hand inequality in (1) is sometimes an equality--
a sphere centered at the origin has hitting probability and
Martin capacity both equal to 1. To see that the constant $1/2$ in~(\ref{eq3.4})
cannot be increased, consider the spherical shell
$$\Lambda_R = \{ x \in \R^d : 1 \leq ||x|| \leq R \} \, .$$
We claim that $\lim_{R \rightarrow \infty} {\rm Cap}_K (\Lambda_R) 
= 2$. Indeed by Proposition \ref{pr1.1}, the Martin capacity of {\em any} 
compact set is at most $2$, while lower bounds tending to 2 for the
capacity of $\Lambda_R$ are established by computing
the energy of the probability measure supported on $\Lambda_R$, with density a constant multiple of
 $ {||x||}^{1-d}$ there.

Next, we pass from the local to the global behavior of Brownian paths. 
Barlow and Taylor (1992) noted that for $d \geq 2$ the set of 
nearest-neighbour lattice points to a Brownian path in $\R^d$
is a subset of $\Z^d$ with dimension 2, using their definition of
dimension which is equivalent to~(\ref{eq2.17}).  This 
is a property of the path near infinity; another such property is 
given by

\begin{pr} \label{pr3.1}
Let $B_d (t)$ denote $d$-dimensional Brownian motion.  Let $\Lambda
\subseteq \R^d$ with $d \geq 3$ and let $\Lambda_1$ be the cubical
fattening of $\Lambda$ defined by
$$\Lambda_1 = \{ x \in \R^d : \exists y \in \Lambda \;~s.t.~\; 
   ||y - x ||_\infty \leq 1 \} .$$
Then a necessary and sufficient condition for the almost sure 
existence of times $t_j \uparrow \infty$ at which $B_d (t_j) \in 
\Lambda_1$ is that ${\rm Cap}_{d-2}^{(\infty)} (\Lambda_1 \cap \Z^d) > 0$.
\end{pr}

The proof is very similar to the proof of Theorem~\ref{th2.1}
and is omitted.

\section{Concluding remarks}

\noindent{\bf 1.} With the exception of Section 5, this paper is concerned with
{\em discrete} Markov chains. Of course the proof of Proposition 
\ref{pr1.1} given in that section extends without difficulty to
some other Markov processes in continuous time, but a classification
of the processes for which this extension is possible is beyond the scope
of this paper. Nevertheless, we do mention explicitly the range
 of a stable subordinator
of index $1/2$, since this range can be viewed as the zero set of
a one-dimensional Brownian motion, and is therefore of wider interest.
\begin{cor} \label{corzeroset}
Let $\{ B(t) \}$ denote standard one-dimensional Brownian 
motion, and let $A$ be any closed set
 in $ (0, \infty) $.  Then
$$
{1 \over 2} {\rm Cap}_K (A) \leq \P [ \exists t \in A \, :\, B(t)=0] 
 \leq {\rm Cap}_K (A) \, \, \mbox{ where}
$$
$$K(s,t) = \left\{ \begin{array}{ll}
               (t/t-s)^{1/2} & \mbox{if $s<t$} \\
               \infty        & \mbox{if $s=t$} \\
                 0           & \mbox{\rm otherwise .}
             \end{array} \right.
$$
\end{cor}
\noindent{\sc Sketch of proof}: Use the obvious estimate
$\P (|B(t)| < \epsilon) \sim 2 \epsilon / \sqrt{2 \pi t}$ as
 $\epsilon \downarrow 0 \, $,
and mimic the proof of Proposition \ref{pr1.1}.

\noindent{\bf 2.} A probabilist might wonder what is gained by capacity
estimates such as Proposition \ref{pr1.1} and Theorem \ref{th2.1},
since the quantity of interest, the hitting probability, is estimated by 
a quantity which appears more complicated. Indeed only in special situations
can the capacity of a set be calculated exactly. Capacity estimates
are useful because of their robustness (see corollaries 2.3 and 2.6,
as well as the proof of the stability of the Nash-Williams recurrence
criterion in Lyons (1992)) and the ease with which they yield
lower bounds for hitting probabilities. Finally, in the continuous setting,
such estimates allow one to exploit the information amassed on capacity by
analysts studying singularities of solutions to PDE's.

\noindent{\bf 3.} The restriction to dimension $d \geq 3$ 
     in Proposition \ref{pr1.1}  is natural since planar
      Brownian motion will hit any measurable set with probability 
        $0$ or $1$. However, by killing the motion at a finite time
       one may obtain a planar version of the proposition.
 
\noindent{\bf 4.} The Martin kernel is most often encountered in constructions
 of the Martin
boundary, where only its asymptotics matter. In Lyons, MacGibbon and Taylor
(1984) the Kernel $G(x,y)/G(0,y)$ is used  to compare the
probability of Brownian motion hitting
a set, to the probability of hitting its projection on a hyperplane.
However, the denominator plays a different role there, as the
Brownian motion is not started at  0, 
and is stopped when it leaves the upper half-space.

\noindent{\bf 5.} The methods of this paper do not seem to yield
upper estimates for the probability that a set will be hit by
the {\em intersection} of the ranges of two Markov chains.
Such estimates were obtained, in a very general setting, in a remarkable
paper by Fitzsimmons and Salisbury (1989). However, the estimates
in that paper required that the initial distribution for each chain be an 
equilibrium measure, so that for fixed initial states only qualitative
(though important) information was obtained. After we showed Tom Salisbury
the statement of Proposition \ref{pr1.1}, he observed that the methods
of his paper with P. Fitzsimmons may be used to estimate
the hitting probability of a set by the intersection of two chains
(with no restrictions imposed on the initial distributions),
in terms of the product of the corresponding Martin kernels.
See Salisbury (1994) for a very readable exposition.

\sc
\noindent Itai Benjamini, Mathematical Sciences Institute, 409 College
Ave., Ithaca, NY, 14853

\noindent Robin Pemantle, Department 
of Mathematics, University of Wisconsin-Madison, Van Vleck Hall, 480 Lincoln
Drive, Madison, WI 53706 . 

\noindent Yuval Peres, Department of
Statistics, 367 Evans Hall University of California, Berkeley, CA 94720

\end{document}